\documentclass[11pt,reqno]{amsart}
\usepackage{amssymb,amsfonts,amscd,amsbsy,color}
\usepackage[mathscr]{eucal}
\usepackage{url}

\newtheorem{theorem}{Theorem}[section]

\theoremstyle{definition}

\newtheorem{conjecture}[theorem]{Conjecture}

\newcommand{\la}{\lambda}

\numberwithin{equation}{section}

\makeatletter
\def\imod#1{\allowbreak\mkern5mu({\operator@font mod}\,\,#1)}
\makeatother

\allowdisplaybreaks

\begin{document}

\title[A positivity conjecture]
{A positivity conjecture related first positive rank and crank moments for overpartitions}

\author{Xinhua Xiong }
\address{Research Institute for Symbolic Computation, Johannes Kepler University\\ A4040 Linz, Austria}
\email{xxiong@risc.uni-linz.ac.at}

\thanks{This research has been by the Austria Science Foundation (FWF) grant SFB F50-06 (Special Research Programm `` Algorithmic and Enumerative Combinatorics'').}

\date{March 20, 2016}

\begin{abstract}
Recently, Andrews, Chan, Kim and Osburn introduced a $q$-series $h(q)$ for the study of the first positive rank and crank moments for overpartitions. They conjectured that for all integers $m \geq 3$,
\begin{equation*}\label{hqcon}
\frac{1}{(q)_{\infty}} (h(q) - m h(q^{m}))
\end{equation*}
has positive power series coefficients for all powers of $q$. Byungchan Kim, Eunmi Kim and Jeehyeon Seo provided a combinatorial  interpretation and proved it is asymptotically true by circle method.  In this note, we show this conjecture is true if $m$ is any positive power of $2$, and we show that in order to prove this conjecture, it is only to prove it for all primes $m$. Moreover we give a stronger conjecture. Our method is very simple and completely different from that of Kim et al.
\end{abstract}

\maketitle

\section{Introduction}
An overpartition \cite{lcc} is a partition in which the first occurrence of each distinct number may be overlined. For example, the $14$ overpartitions of $4$ are
\begin{equation*}
\begin{gathered}
4, \overline{4}, 3+1, \overline{3} + 1, 3 + \overline{1},
\overline{3} + \overline{1}, 2+2,\\ \overline{2}
+ 2, 2+1+1, \overline{2} + 1 + 1, 2+ \overline{1} + 1, \\
\overline{2} + \overline{1} + 1, 1+1+1+1, \overline{1} + 1 + 1 +1.
\end{gathered}
\end{equation*}

Let $\overline{N}(n,m)$ denote the number of overpartitions of $n$ whose rank is $m$ and $\overline{M}(n,m)$ the number of overpartitions of $n$ whose first residual crank is $m$. Andrews, Chan, Kim and Osburn \cite{acko} defined the positive rank and crank moments for overpatitions: 

\begin{equation*}
\overline{N}_{k}^{+}(n) := \sum_{m=1}^{\infty} m^{k} \overline{N}(m,n)
\end{equation*}

\noindent and

\begin{equation*}
\overline{M}_{k}^{+}(n) := \sum_{m=1}^{\infty} m^{k} \overline{M}(m,n).
\end{equation*}

\noindent They proved the inequality
\begin{equation} \label{oddrcover}
\overline{M}_{1}^{+}(n) > \overline{N}_{1}^{+}(n)
\end{equation}
holds for all $n \geq 1$. They also gave the difference 
$ \overline{M}_{1}^{+}(n) - \overline{N}_{1}^{+}(n)$ a combinatorial interpretation. In order to prove (\ref{oddrcover}),  Andrews, Chan, Kim and Osburn \cite{acko} introduced the function
\begin{equation*} \label{hqdefn}
 h(q) := \sum_{n=1}^{\infty} \frac{(-1)^{n+1} q^{n(n+1)/2}}{1-q^n}, 
 \end{equation*}
 and conjectured that
$$\label{con}
 \frac{1}{(q)_{\infty}} ( h(q) - m h(q^m) )
$$
 has positive coefficients for all $m\geq 3$, where we use the standard $q$-series notation, $(a)_{\infty} = (a;q)_{\infty} = \prod_{n=1}^{\infty} (1-aq^{n-1} ).$  In \cite{kks}, Byungchan Kim, Eunmi Kim and Jeehyeon Seo provided a combinatorial interpretation for the coefficients of $ \frac{1}{(q)_{\infty}} h( q^m )$. According to their definition, a $m$-string in an ordinary partition is the parts consisting of $m(1+k)$, $m(3+k)$,  $\ldots$, $m(2j-1+k)$ with a positive integer $j$ and a nonnegative integer $k \le j$, and a weight of $m$-string is 1 if $k=0$ or $j$, and 2, otherwise. They defined $C_m (n)$ as the weighted sum of the number of $m$-strings along the partitions of $n$, i.e.
$$
 C_m (n) = \sum_{\la \vdash n} \sum_{\substack{\text{$\pi$ is a} \\ \text{$m$-string of $\la$}}} \operatorname{wt} (\pi).
$$
It is clear that
$$
 \frac{1}{(q)_{\infty}} h(q^m) = \sum_{n \ge 1} C_m (n) q^n.
$$
So the conjecture of Andrews, Chan, Kim and Osburn \cite{acko} can be intepretated as there are more (weighted count of) $1$-strings than $m$ times of (weighted count of) $m$-strings along the partitions of $n$. In the paper \cite{kks}, Byungchan Kim, Eunmi Kim and Jeehyeon Seo, by using the circle method of Wright \cite{w} and some results
from \cite{bm2},  proved the conjecture of Andrews-Chan-Kim-Osburn is asymptotically true. In this note, we will prove this conjecture is true if $m$ is any power of $2$. Moreover, we show that in order to prove the conjecture, it is only to prove it is true for all primes $m$.  Here are our main results.
\begin{theorem}\label{main}
For all integer $m\geq 2$,
$$
\frac{1}{(q)_{\infty}}(h(q)-2^{m}h(q^{2^{m}})
$$
has positive power series coefficients for all positive powers of $q$.
\end{theorem}

\begin{theorem}\label{second main}
Suppose for a prime $p$, 
$$
\frac{1}{(q)_{\infty}}(h(q)-ph(q^{p})
$$
has positive power series coefficients for all positive powers of $q$.
Then for all integer $m\geq 2$,
$$
\frac{1}{(q)_{\infty}}(h(q)-p^{m}h(q^{p^{m}})
$$
has positive power series coefficients for all positive powers of $q$.
\end{theorem}
We don't know why primes appeared here, in Section 3, we give a stronger conjecture related primes.

\section{ Proof of Theorem 1.1. and Theorem 1.2}
\noindent For a prime $p\geq 1$, we define the function
$
M_{p}(q)=\frac{1}{(q)_{\infty}}(h(q)-ph(q^p)),
$
\noindent then we have 
$M_{2}(q)$ has positive power series coefficients for all positive powers of $q^{n}$ for all $n$ except that the coefficients of $q^{2}$ and $q^{4}$ are zero. This lemma is the Corollary 2.4 of \cite{acko}.

{\it Proof of Theorem 1.1.} $m\geq 2$,
\begin{multline*} \label{1}
\begin{aligned}
&M_{2^{m}}(q)
 = \frac{1}{(q)_{\infty}}(h(q)-2^{m}h(q^{2^{m}})) \\
& = \frac{1}{(q)_{\infty}}(h(q)-2h(q^2)+2h(q^{2})-4h(q^4)+\cdot\cdot\cdot
+2^{m-1}h(q^{2^{m-1}})-2^{m}h(q^{2^{m}}))\\
& =\frac{h(q)-2h(q^2)}{(q)_{\infty}}+2\frac{h(q^{2})-2h(q^4)}{(q)_{\infty}}+\cdot\cdot\cdot
+2^{m-1}\frac{h(q^{2^{m-1}})-2h(q^{2^{m}})}{(q)_{\infty}}\\
&=\frac{h(q)-2h(q^2)}{(q)_{\infty}}+
\frac{2}{(1-q)(1-q^{3})(1-q^{5})\cdot\cdot\cdot}\cdot\frac{h(q^{2})-2h(q^{4})}{\prod^{\infty}_{n=1}(1-q^{2n})}+\cdot\cdot\cdot \\
&+\frac{2^{m-1}}{\prod_{n\neq 2^{m-1}k, k\geq 1}(1-q^{n})}\cdot \frac{h(q^{2^{m-1}})-2h(q^{2^{m}})}{\prod_{n=1}(1-q^{2^{m-1}n})}\\
&=M_{2}(q)+\frac{2}{(1-q)(1-q^{3})(1-q^{5})\cdot\cdot\cdot}\cdot M_{2}(q^{2}) +  \cdot\cdot\cdot 
+ \frac{2^{m-1}}{\prod_{n\neq 2^{m-1}k, k\geq 1}(1-q^{n})}\cdot M_{2}(q^{2^{m-1}}).
\end{aligned}
\end{multline*}
\noindent By the lemma above, $M_{2}(q)$ has positive coefficients of $q^{n}$ for all $n$ with except 2 and 4, therefore for all $m\geq 2$,  $M_{2}(q^{2^{m-1}})$ has positive coefficients except that the coefficients of $q^{2^{m+1}}$ and $q^{2^{m+2}}$, but the sum of the left hand will have positive coefficients of $q^{n}$ for all $\geq 1$.
\qed

{\it Proof of Theorem 1.2.}
For $p\geq 3$ a prime and $m\geq 2$, 
\begin{equation*} 
\begin{aligned}
&M_{p^{m}}(q)
 = \frac{1}{(q)_{\infty}}(h(q)-p^{m}h(q^{p^{m}})) \\
& = \frac{1}{(q)_{\infty}}(h(q)-ph(q^p)+ph(q^{p})-p^{2}h(q^{p^{2}})+\cdot\cdot\cdot
+p^{m-1}h(q^{p^{m-1}})-p^{m}h(q^{p^{m}}))\\
& =\frac{h(q)-ph(q^p)}{(q)_{\infty}}+p\frac{h(q^{p})-ph(q^{p^{2}})}{(q)_{\infty}}+\cdot\cdot\cdot
+p^{m-1}\frac{h(q^{p^{m-1}})-ph(q^{p^{m}})}{(q)_{\infty}}\\
&=\frac{h(q)-ph(q^p)}{(q)_{\infty}}+
\frac{p}{\prod_{n\neq pk, k\geq 1}(1-q^{n})}\cdot\frac{h(q^{p})-2h(q^{p^{2}})}{\prod^{\infty}_{n=1}(1-q^{pn})}+\cdot\cdot\cdot \\
&+\frac{p^{m-1}}{\prod_{n\neq p^{m-1}k, k\geq 1}(1-q^{n})}\cdot \frac{h(q^{p^{m-1}})-ph(q^{p^{m}})}{\prod_{n=1}(1-q^{p^{m-1}n})}\\
&=M_{p}(q)+\frac{p}{\prod_{n\neq pk, k\geq 1}(1-q^{n})}\cdot M_{p}(q^{p}) +  \cdot\cdot\cdot 
+ \frac{p^{m-1}}{\prod_{n\neq p^{m-1}k, k\geq 1}(1-q^{n})}\cdot M_{p}(q^{p^{m-1}}).
\end{aligned}
\end{equation*}
Since each summand of the right hand side of the above has positive coefficients of $q^{n}$ for all positive integers $n$, $M_{p^{m}}$ will have positive coefficients
of $q^{n}$ for all integers $n$.
\qed

{\it Remark.}
By our method, it can be easily seen that if the conjecture is true for all primes $m=p$, then it is true for any other natural numbers $m$. For example, consider the case $m=6$, 
\begin{equation*} \label{1}
\begin{aligned}
&\frac{h(q)-6h(q^{6})}{(q)_{\infty}}\\
& = \frac{h(q)-2h(q^{2})}{(q)_{\infty}}+\frac{2h(q)-6h(q^{6}))}{(q)_{\infty}}\\
&=M_{2}(q)+\frac{2}{(1-q)(1-q^{2})(1-q^{3}\cdot\cdot\cdot)}\cdot \\
&\frac{h(q^{2})-3h(q^{6})}{(1-q^{2})(1-q^{4})(1-q^{6})\cdot\cdot\cdot}\\
&=M_{2}(q) + \frac{2}{(1-q)(1-q^{3})(1-q^{5})\cdot\cdot\cdot}\cdot M_{3}(q^{2}).
\end{aligned}
\end{equation*}
\noindent We see that the positivity of coefficients of the power series $M_{2}(q)$ and $M_{3}(q)$ will imply the positivity of the coefficients of the power series of $M_{6}(q)$. 

\section{A Stronger Conjecture }
Kim et al. proved that the conjecture of Andrews et al is asymptotically true by using circle method, But Andrews et al. originally expected to find $q$-theoretic or combinatorial proofs for this conjecture.  Here basing on the numerical results, we make the following stronger conjecture. We also expected find $q$-theoretic or combinatorial proofs for this conjecture. We can easily see that this conjecture implies the conjecture of Andrews et al.
\begin{conjecture}
Let $p_{1}>p_{2} \geq2$ be two primes.  Then the function
$$
 \frac{1}{(q)_{\infty}} ( p_{1}h(q^{p_{1}}) - p_{2}h(q^{p_{2}}) )$$
has positive power series coefficients of $q^{n}$ for all $n\geq p_{2}$ and has nonnegative power series coefficients of $q^{n}$ for all $n\geq 1$.
\end{conjecture}
We verified this conjecture for the first $100,000$ coefficients of the power series for each prime pair cases which are less than $50$ by using Mathematica. We provide some coefficients of power series of  $\frac{1}{(q)_{\infty}} ( h(q) - m h(q^m) )$ for small prime $m$, which are also obtained by using Mathematica.
\begin{equation*}
\begin{aligned}
 &\frac{1}{(q)_{\infty}} ( h(q) - 3 h(q^3) )\\
 &=q + 2 q^2 + 3 q^4 + 3 q^5 + 4 q^6 + 5 q^7 + 9 q^8 + 10 q^9 \\
 &+16 q^{10} + 19 q^{11} + 26 q^{12} + 33 q^{13} + 46 q^{14} + 56 q^{15} +\cdot\cdot\cdot \\
 \end{aligned}
 \end{equation*} 
  \begin{equation*}
\begin{aligned} & \frac{1}{(q)_{\infty}} ( h(q) - 5 h(q^5) \\
& =q + 2 q^2 + 3 q^3 + 6 q^4 + 4 q^5 + 11 q^6 + 13 q^7 + 21 q^8 \\
&+27 q^9 + 36 q^{10} + 46 q^{11} + 67 q^{12} + 82 q^{13} + 111 q^{14} \\
&+ 141 q^{15} +\cdot\cdot\cdot
 \end{aligned}
 \end{equation*} 
  \begin{equation*}
\begin{aligned} & \frac{1}{(q)_{\infty}} ( h(q) - 7 h(q^7) )=q + 2 q^2 + 3 q^3 + 6 q^4 + 9 q^5 \\
&+ 16 q^6 + 16 q^7 + 29 q^8 + 
 38 q^9 + 55 q^{10} \\
 &+ 71 q^{11} + 103 q^{12 }+ 130 q^{13} + 174 q^{14 }+ 
 225 q^{15} +\cdot\cdot\cdot\\
  \end{aligned}
 \end{equation*}

\begin{equation*}
\begin{aligned}& \frac{1}{(q)_{\infty}} ( h(q) - 11 h(q^{11}) ) =q + 2 q^2 + 3 q^3 + 6 q^4 + 9 q^5\\
& +16q^6 + 23q^7 + 36 q^8 
 +52q^9 + 76q^{10 }+ 95q^{11}\\
 & + 141 q^{12} + 185 q^{13} +253 q^{14} + 331 q^{15}+\cdot\cdot\cdot \\
  \end{aligned}
 \end{equation*} 

 \begin{equation*}
\begin{aligned}
 & \frac{1}{(q)_{\infty}} ( h(q) - 13 h(q^{13}) )=q + 2 q^2 + 3 q^3 + 6 q^4 + 9 q^5\\
 & + 16 q^6 + 23 q^7 + 36 q^8 + 
 52 q^9 + 76 q^{10}\\
 & + 106 q^{11} + 152 q^{12} + 192 q^{13} + 273 q^{14} +360q^{15}+ \cdot\cdot\cdot\\
  \end{aligned}
 \end{equation*} 

 \begin{equation*}
\begin{aligned} & \frac{1}{(q)_{\infty}} ( h(q) - 17 h(q^{17}) )=q + 2 q^2 + 3 q^3 + 6 q^4 + 9 q^5\\
 &+16q^6 + 23 q^7 + 36 q^8 + 
52 q^9 + 76q^{10}\\
 & + 106 q^{11} + 152 q^{12} + 207 q^{13} + 286 q^{14 }+ 
 386 q^{15} + \cdot\cdot\cdot\\
  \end{aligned}
 \end{equation*}

\section*{Acknowledgments}
The author would like to thank Professor Peter Paule  for his valuable comments on an earlier version of this paper and encouragements.
The author also would like to thank Professor Gorge Andrews  for his suggestion of the rewriting some parts of this note.

\end{document}